\documentclass[11pt,leqno]{article}

\usepackage{amsmath,amsfonts,amscd,amssymb,theorem}

\long\def\comment#1\endcomment{}

\comment
\endcomment

\makeatletter
\begingroup
\gdef\th@dotted{\normalfont\itshape
  \def\@begintheorem##1##2{%
        \item[\hskip\labelsep \theorem@headerfont ##1\ ##2.]}%
\def\@opargbegintheorem##1##2##3{%
   \item[\hskip\labelsep \theorem@headerfont ##1\ ##2\ (##3).]}}
\endgroup
\makeatother

\theoremstyle{dotted}

\newtheorem{theorem}{Theorem}[section]
\newtheorem{lemma}[theorem]{Lemma}

\newtheorem{prop}[theorem]{Proposition}
\newtheorem{corr}[theorem]{Corollary}

\makeatletter
\begingroup
\gdef\th@upshape{\normalfont
  \def\@begintheorem##1##2{%
        \item[\hskip\labelsep \theorem@headerfont ##1\ ##2.]}%
\def\@opargbegintheorem##1##2##3{%
   \item[\hskip\labelsep \theorem@headerfont ##1\ ##2\ (##3).]}}
\endgroup
\makeatother

\theoremstyle{upshape}

\newtheorem{defn}[theorem]{Definition}
\newtheorem{remark}[theorem]{Remark}
\newtheorem{exa}[theorem]{Example}

\makeatletter
\renewcommand{\subsection}{\@startsection{subsection}{2}{0pt}{-3ex
plus -1ex minus -0.2ex}{-2mm plus -0pt minus
-2pt}{\normalfont\bfseries}} 
\renewcommand{\subsubsection}{\@startsection{subsubsection}{3}{0pt}{-3ex
plus -1ex minus -0.2ex}{-2mm plus -0pt minus
-2pt}{\normalfont\bfseries}} 
\makeatother

\makeatletter
\@addtoreset{equation}{section}
\makeatother

\newcommand{\cntrct}                % contraction with a vector field
{\hspace{2pt}\raisebox{1pt}{\text{$\lrcorner$}}\hspace{2pt}}

\newcommand{\proof}[1][Proof.]{\smallskip\noindent{\em #1}}
\def\endproof{\hfill\ensuremath{\square}\par\medskip}

\def\eqref#1{\thetag{\ref{#1}}}

\let\latexref=\ref
\def\ref#1{{\normalfont{\latexref{#1}}}}

\newcommand{\wt}{\widetilde}
\newcommand{\wh}{\widehat}

\newcommand{\idot}{{\:\raisebox{1pt}{\text{\circle*{1.5}}}}}
\newcommand{\hdot}{{\:\raisebox{3pt}{\text{\circle*{1.5}}}}}
\newcommand{\eps}{\varepsilon}
\renewcommand{\phi}{\varphi}

\newcommand{\Spec}{\operatorname{Spec}}
\newcommand{\Fun}{\operatorname{Fun}}

\newcommand{\Ext}{\operatorname{Ext}}
\newcommand{\Tor}{\operatorname{Tor}}
\newcommand{\Coker}{\operatorname{Coker}}
\newcommand{\Ker}{\operatorname{Ker}}

\newcommand{\cchar}{\operatorname{\sf char}}
\newcommand{\ppt}{{\sf pt}}
\newcommand{\id}{\operatorname{\sf id}}

\newcommand{\D}{{\cal D}}

\newcommand{\R}{\mathbb{R}}
\newcommand{\Z}{\mathbb{Z}}
\newcommand{\Q}{\mathbb{Q}}
\newcommand{\N}{\mathbb{N}}

\newcommand{\K}{\mathbb{K}}

\newcommand{\lotimes}{\overset{\sf\scriptscriptstyle L}{\otimes}}

\newcommand{\C}{\mathcal{C}}

\newcommand{\tr}{\operatorname{\sf tr}}
\newcommand{\e}{\operatorname{\sf e}}

\newcommand{\hash}{\sharp}

\newcommand{\I}{{\sf I}}

\newcommand{\vH}{\check{H}}

\title{Cartier isomorphism for unital associative algebras}

\author{D. Kaledin\thanks{Partially supported by RScF, grant
    14-50-00005, and the Dynasty Foundation award}}

\begin{document}

\maketitle

\tableofcontents

\section*{Introduction.}

For a smooth affine scheme $X = \Spec A$ of finite type over a
perfect field $k$ of positive characteristic, the {\em Cartier
  isomorphism} identifies the de Rham cohomology of $X$ with the
spaces of differential forms. More precisely, for any $i \geq 0$,
let $\Omega^\hdot(X)$ be space of global $i$-forms on $X$, and let
$$
B^\hdot(X),Z^\hdot(X) \subset \Omega^\hdot(X)
$$
be the subspaces of exact resp. closed forms. Then there exists a
canonical isomorphism
$$
C:Z^\hdot(X)/B^\hdot(X) \cong \Omega^\hdot(X).
$$
This isomorphism is Frobenius-semilinear with respect to the
Frobenius endomorphism of the field $k$. It is functorial and
compatible with localizations, so that it induces an isomorphism of
Zariski sheaves on $X$ (and in this from, it generalizes to
non-affine schemes).

The classic theorem of Hochschild, Kostant and Rosenberg identifies
the spaces $\Omega^\hdot(X)$ with the Hochschild homology groups
$HH_\idot(A)$ of the algebra $A$ -- for any $i \geq 0$, we have a
canonical isomorphism
$$
\Omega^i(X) \cong HH_i(A).
$$
The de Rham differential $d$ coincides under this identification
with the Connes-Tsygan differential $B:HH_\idot(A) \to HH_{\idot +
  1}(A)$. However, both the Hochschild homology groups and the
Connes-Tsygan differential exist in a larger generality -- they are
defined for an arbitrary unital associative algebra $A$. It is
natural to ask whether the Cartier isomorphism exists in this larger
generality.

\medskip

The goal of this paper is to give a positive answer to this
question, provided the characteristic of the base field is not
$2$. Namely, for any associative unital $k$-algebra $A$, we
construct functorial groups $BHH_\idot(A)$, $ZHH_\idot(A)$ and maps
$$
\begin{CD}
BHH_\idot(A) @>{\xi}>> ZHH_\idot(A) @>{\zeta}>> HH_\idot(A)
@>{\beta}>> BHH_{\idot+1}(A)
\end{CD}
$$
such that that Connes-Tsygan differential factors as $B = \zeta
\circ \xi \circ \beta$, and the map $\xi$ fits into a natural
long exact sequence
$$
\begin{CD}
BHH_\idot(A) @>{\xi}>> ZHH_\idot(A) @>{C}>> HH_\idot(A)
@>>>,
\end{CD}
$$
with a certain natural Frobenus-semilinear map $C$. This map $C$ is
our non-commutative generalization of the Cartier isomorphism. We
then prove that in the Hochschild-Kostant-Rosenberg case -- that is,
$A$ is commutative of finite type and $X = \Spec A$ is smooth -- our
generalized objects reduce to the classic ones: $BHH_\idot(A)$ and
$ZHH_\idot(A)$ become $B^\hdot(X)$ and $Z^\hdot(X)$, the maps
$\alpha$ and $\beta$ are the natural embeddings, and $C$ is the
usual Cartier isomorphism.

Note that while the Cartier isomorphisms of the individual
cohomology groups are the most one can get in the general case,
stronger results can be obtained if one imposes retrictions on the
variety $X$. For example, if one assumes that $X$ is liftable to the
ring $W_2(k)$ of second Witt vectors of $k$, and that $\dim X$ is
less than $\cchar k$, then $C$ can be lifted to an isomorphism
$$
F_*\Omega^\hdot(X) \cong \bigoplus_i\Omega^i(X^{(1)})[-i]
$$
in the derived category of coherent sheaves on the Frobenius twist
$X^{(1)}$, where $F:X \to X^{(1)}$ is the Frobenius map. This forms
the basis of the famous paper \cite{DI} of P. Deligne and
L. Illusie.  A non-commutative generalization of the Deligne-Illusie
isomorphism has been already constructed in \cite{ka0} (under the
assumptions that $A$ lifts to $W_2(k)$, and the homological
dimension of the category of $A$-bimodules is less than $2\cchar
k$). Therefore in this paper, we concentrate on the results that
require no assumptions on the associative algebra $A$.

Another topic we avoid in this paper is the so-called {\em conjugate
  spectral sequence} -- namely, the spectral sequence generated by
the canonical filtration on the de Rham complex
$\Omega^\hdot(X)$. In the classic case, this spectral sequence
converges to de Rham cohomology of $X$, and the Cartier isomorphism
identifies the first term of the spectral sequence with the
Frobenius twist of the Hodge cohomology of $X$. The corresponding
statement for cyclic homology is also true, in a sense, but it is
significantly more delicate because of possible convergence issues
-- in effect, what the conjugate spectral sequence converges to is
not the usual periodic cyclic homology but a different homology
theory. We note that already in the commutative case, a similar
phenomenon occurs when one studies singular schemes and derived de
Rham complexes, see e.g. the recent work of A. Beilinson \cite{Be}
and B. Bhatt \cite{Bh}. We will return to this elsewhere.

\medskip

The paper is organized as follows. In Section~\ref{pre.sec}, we
recall the necessary generalities on cyclic categories, cyclic
homology and so on. The material is pretty standard. We follow the
exposition given previously in \cite{ka0} and use more-or-less the
same notation and conventions. In Section~\ref{alg.sec}, we recall
how cyclic homology works for associative algebras -- again, this
material is very standard. Section~\ref{tra.sec} is concerned with a
slightly non-standard part of the cyclic homology story, so we give
complete proofs. Having finished with the preliminaries, we can turn
to our subject. In Section~\ref{car.sec}, we define our
non-commutative Cartier map $C$ and the whole package that comes
with it -- the groups $BHH_\idot(A)$, $ZHH_\idot(A)$ and all the
natural maps between them. The map $C$ appears at the very end of
the Section in Definition~\ref{car.def}. Finally, in
Section~\ref{comp.sec}, we prove the comparison result stating that
in the HKR case, our generalized Cartier maps reduces to the
classical one (this is Proposition~\ref{car.comp}).

\subsection*{Acknowledgements.} I am grateful to the referee for
useful comments and suggestions (in particular, Remark~\ref{Efi}).

I also feel that in present circumstances, it is my duty to
acknowledge my gratitude to the Dynasty Foundation, for the
financial support extended over the years both to me and to many
other mathematicians in Russia. As far as I know, until its forcible
closure by the powers-that-be in June, 2015, Dynasty Foundation had
been the only source of funding for mathematics in Russia that was
not related to the government of the Russian Federation, and
coincidentally, the only source of funding where one could be sure
that the money was free from any criminal taint.

\section{Preliminaries.}\label{pre.sec}

Our approach to cyclic homology is based on homology of small
categories, as in \cite{C} and \cite{FT}. Our notation is as
follows. For any small category $I$ and ring $R$, we denote by
$\Fun(I,R)$ the category of functors from $I$ to the category of
left $R$-modules. This is an abelian category. We denote its derived
category by $\D(I,R)$. If $I = \ppt$ is the point category, then
$\Fun(\ppt,R)$ is the category of left $R$-modules, and $\D(\ppt,R)$
is its derived category $\D(R)$. For any functor $\gamma:I' \to I$
between small categories, we denote by $\gamma^*:\Fun(I,R) \to
\Fun(I',R)$ the natural pullback functor, and we denote by
$\gamma_!,\gamma_*:\Fun(I',R) \to \Fun(I,R)$ its left and
right-adjoint, namely, the left and right Kan extension along
$\gamma$. The derived functors
$L^\hdot\gamma_!,R^\hdot\gamma_*:\D(I',R) \to \D(I,R)$ are left
resp. right-adjoint to the pullback functor $\gamma^*:\D(I,R) \to
\D(I',R)$. The {\em homology complex} of a small category $I$ with
coefficients in an object $E \in \D(I,R)$ is then given by
$$
C_\idot(I,E) = L^\hdot\tau_!E \in \D(R),
$$
where $\tau:I \to \ppt$ is the projection to the point. The homology
groups $H_\idot(I,E)$ are the homology groups of the complex
$C_\idot(I,R)$.

The small category one needs for applications to cyclic homology is
the {\em cyclic category} $\Lambda$ introduced by A. Connes
\cite{C}. It has many equivalent definitions, for example the
following:
\begin{itemize}
\item objects of $\Lambda$ are finite cellular decompositions of the
  circle $S^1$, with $[n] \in \Lambda$ standing for the
  decomposition with $n$ cells of dimension $0$,
\item morphisms from $[n]$ to $[m]$ are homotopy classes of
  continuous cellular maps $f:S^1 \to S^1$ of degree $1$ such that
  the universal covering map $\wt{f}:\R \to \R$ is order-preserving.
\end{itemize}
For any $[n] \in \Lambda$, $0$-cells in the corresponding cellular
decomposition are called {\em vertices}, and the set of vertices is
denoted $V([n])$. For alternative definitions of $\Lambda$, see
e.g. \cite{Lo}.

We will also need the following modification of the category
$\Lambda$. Fix an integer $n \geq 1$, and consider the $n$-fold
covering map $S^1 \to S^1$. Then a cellular decomposition of $S^1$
corresponding to an object $[m] \in \Lambda$ induces a cellular
decomposition of its $n$-fold cover, thus gives an object $[nm] \in
\Lambda$, and the deck transformation gives a map $\sigma:[nm] \to
[nm]$ in $\Lambda$ such that $\sigma^n=\id$. The {\em category
  $\Lambda_n$} is the category of such objects $[nm] \in \Lambda$
and morphisms between them that commute with $\sigma$. Forgetting
$\sigma$ then gives a functor
$$
i_n:\Lambda_n \to \Lambda,
$$
and on the other hand, taking the quotient by $\sigma$ gives a
functor
\begin{equation}\label{pi.n}
\pi^n:\Lambda_n \to \Lambda.
\end{equation}
Objects of $\Lambda_n$ are numbered by integers $m \geq 1$; we
denote the object corresponding to $m$ by $[m] \in \Lambda$, so that
$\pi^n([m])=[m]$ and $i_n([m])=[mn]$. We note that the functor
$\pi^n$ of \eqref{pi.n} is a bifibration in the sense of \cite{SGA}.
Each fiber of this bifibration is the groupoid $\ppt_n$ with one
object with automorphism group $\Z/n\Z$ generated by $\sigma$. Of
course, for $n=1$, we have $\Lambda_1=\Lambda$ and $i_1 = \pi^1 =
\id$.

As usual, we denote by $\Delta$ the category of finite non-empty
totally ordered sets, with $[n] \in \Delta$ being the set with $n$
elements, and we denote by $\Delta^o$ the opposite category. One
shows -- see e.g. \cite[Subsection 1.5]{ka0} that uses exactly the
same notation as this paper -- that there exists a natural functor
$$
j:\Delta^o \to \Lambda
$$
sending $[n] \in \Delta^o$ to $[n] \in \Lambda$ (in effect,
$\Delta^o$ is equivalent to the category $[1] \setminus \Lambda$ of
objects $[n] \in \Lambda$ equipped with a map $[1] \to [n]$, and $j$
is the functor that forgets the map). For any $n \geq 1$, we also
have an analogous functor $j_n:\Delta^o \to \Lambda_n$, and $\pi^n
\circ j_n \cong j$.

\begin{defn}
For any $n \geq 1$, any ring $R$, and any object $E \in
\D(\Lambda_n,R)$, the {\em Hochschild homology} $HH_\idot(E)$
resp. the {\em cyclic homology} $HC_\idot(E)$ are given by
$$
HH_\idot=H_\idot(\Delta^o,j_n^*E), \qquad HC_\idot(E) =
H_\idot(\Lambda_n,E).
$$
\end{defn}

For any $E \in \Fun(\Delta^o,R)$, the homology $H_\idot(\Delta^o,E)$
can be computed by the standard complex $CH_\idot(E)$ of the
simplicial $R$-module $E$, namely, $CH_i(E) = E([i+1])$, $i \geq 0$,
with the differential given by the alternating sum of the face
maps. By abuse of notation, for any $E \in \Fun(\Lambda_n,R)$, we will
denote $CH_\idot(E) = CH_\idot(j^*_nE)$, and call it the {\em
  Hochschild homology complex} of the object $E$.

To analyze cyclic homology, it is convenient to do the
following. For any $[n] \in \Lambda$, denote by $K_\idot([n])$ the
cellular chain complex of the circle $S^1$ computing its homology
with coefficients in $\Z$. This is functorial with respect to
cellular maps, and the homology does not depend on the choice of a
cellular decomposition, so that we obtain an exact sequence
\begin{equation}\label{4.term}
\begin{CD}
\Z @>{\kappa_1}>> \K_1 @>>> \K_0 @>{\kappa_0}>> \Z
\end{CD}
\end{equation}
in $\Fun(\Lambda,\Z)$, where $\Z$ stands for the constant
functor. For any $n \geq 1$, we can apply the pullback functor
$i_n^*$ and obtain an analogous exact sequence in
$\Fun(\Lambda_n,\Z)$. For any $E \in \Fun(\Lambda_n,R)$, denote
$$
\K_i(E) = i_n^*\K_i \otimes E, \qquad i=0,1.
$$
Then $\K_\idot(E)$ is a complex in $\Fun(\Lambda_n,R)$ of length
$2$, with homology objects in degrees $0$ and $1$ canonically
identified with $E$.

\begin{lemma}\label{la.ho}
\begin{enumerate}
\item For any $n \geq 1$ and $E \in D(\Lambda,R)$, the adjunction
  maps
$$
h_n:HH_\idot(i_n^*E) \to HH_\idot(E), \qquad c_n:HC_\idot(i_n^*E) \to
HC_\idot(E)
$$
are isomorphisms.
\item For any $n \geq 1$ and $E \in \Fun(\Lambda_n,R)$, we have
  $HH_\idot(\K_1(E)) = 0$, and we have a natural isomorphism
  $HH_\idot(E) \cong HC_\idot(\K_0(E))$, so that we have an
  identification
\begin{equation}\label{hh.hc}
HH_\idot(E) \cong HC_\idot(\K_\idot(E)).
\end{equation}
\end{enumerate}
\end{lemma}

\proof{} The isomorphism $h_n$ in \thetag{i} is the so-called {\em
  edgewise subdivision} isomorphism that goes back to Quillen and
Segal; for an independent proof of \thetag{i} with exactly the same
notation as here, see \cite[Lemma 1.14]{ka0}. For \thetag{ii}, see
\cite[Subsection 1.6]{ka0} and specifically \cite[Lemma 1.10]{ka0}.
\endproof

Now for any $E \in \Fun(\Lambda_n,R)$, consider the natural maps
\begin{equation}\label{kappa}
\kappa_0:\K_0(E) \longrightarrow E, \qquad \kappa_1:E
\longrightarrow \K_1(E)
\end{equation}
induced by the maps \eqref{4.term}, and let
\begin{equation}\label{B.0}
B = \kappa_1 \circ \kappa_0:\K_0(E) \to E \to \K_1(E)
\end{equation}
be their composition. We can form a natural bicomplex
\begin{equation}\label{reso}
\begin{CD}
@.@.@.\K_0(E)\\
@.@.@.@AAA\\
@.@.\K_0(E) @>{B}>> \K_1(E)\\
@.@.@AAA\\
@. \K_0(E) @>{B}>> \K_1(E)\\
@.@AAA\\
@>{B}>> \K_1(E)
\end{CD}
\end{equation}
The total complex of this bicomplex is a resolution of the object
$E$. If we denote by $u$ the endomorphism of \eqref{reso} obtain by
shifting to the left and downward by $1$ term, then by virtue of
\eqref{hh.hc}, one of the two standard spectral sequences associated
to a bicomplex reads as
\begin{equation}\label{sp.se}
HH_\idot(E)[u] \Rightarrow HC_\idot(\Lambda_l,E),
\end{equation}
where the left-hand side is shorthand for ``formal polynomials in
one variable $u$ of homological degree $2$ with coefficients in
$HH_\idot(E)$''. This is known as the {\em
Hochschild-to-cyclic}, or {\em Hodge-to-de Rham} spectral sequence.
The first non-trivial differential
\begin{equation}\label{c.ts}
B:HH_\idot(E) \to HH_{\idot+1}(E)
\end{equation}
is known as the {\em Connes-Tsygan differential}, or sometimes {\em
  Rinehart differential}. In terms of the identification
\eqref{hh.hc}, this differential is induced by the map
\begin{equation}\label{B.gen}
B:\K_\idot(E) \to \K_\idot(E)[-1]
\end{equation}
which is in turn induced by the natural map $B$ of \eqref{B.0} (this
is why we use the same notation for all three). Moreover, $B$ can be
lifted to a map of complexes
\begin{equation}\label{B.ch}
B:CH_\idot(E) \to CH_\idot(E)[-1]
\end{equation}
so that $B^2=0$ (see e.g. \cite{Lo}).

If the ring $R$ is commutative, then the categories
$\Fun(\Lambda,R)$, $\Fun(\Delta^o,R)$ acquire natural pointwise
tensor products, and their derived categories in turn acquite the
derived tensor product $\lotimes$. It is well-known that the
Hoch\-schild homology functor is multiplicative: for any two objects
$E,E' \in \D(\Lambda,R)$, we have a natural K\"unneth
quasiisomorphism
\begin{equation}\label{kunn}
CH_\idot(E) \lotimes CH_\idot(E') \cong CH_\idot(E \lotimes E'),
\end{equation}
where by abuse of notation we use $CH_\idot(-)$ to denote the
corresponding object in the derived category $\D(R)$.

\begin{lemma}\label{B.kun}
The Connes-Tsygan differential $B$ of \eqref{c.ts} is compatible
with the K\"unneth quasiisomorphism \eqref{kunn} -- that is, we have
$$
B_{E \lotimes E'} = B_{E} \otimes \id + \id \otimes B_{E'}.
$$
\end{lemma}

This is very well-known if $R$ contains $\Q$. However, we will need
the general case. Since we could not find a convenient reference, we
give a proof.

\begin{defn}
A {\em mixed complex} in an abelian category $\C$ is a pair $\langle
\K_\idot,B \rangle$ of a complex $\K_\idot$ in $\C$ and a map
$B:\K_\idot \to \K_\idot[-1]$ such that $B^2=0$. A map of mixed
complexes is a {\em quasiisomorphism} if it is a quasiismorphism of
the underlying complexes in $\C$.
\end{defn}

\begin{exa}
For any ring $R$ and any $E \in \Fun(\Lambda,R)$, the complex
$\K_\idot(E)$ with the differential $B$ of \eqref{B.gen} is a mixed
complex in $\Fun(\Lambda,R)$.
\end{exa}

Inverting quasiisomorphisms of mixed complexes, we obtain a
triangulated category $\D_{mix}(\C)$. If $\C$ is a tensor category,
then one defines the tensor product of two mixed complexes by
$$
\langle \K_\idot,B \rangle \otimes \langle K'_\idot,B' \rangle =
\langle \K_\idot \otimes K'_\idot, B \otimes \id + \id \otimes B'
\rangle,
$$
and the derived functor of this tensor product obviously turns
$\D_{mix}(\C)$ into a tensor triangulated category (at least when
every complex in $\C$ has a homotopically flat replacement, and this
is the only case that we will need). If $\C'$ is a module category
over $\C$, then the same formula defines an external tensor product
$\langle \K_\idot,B \rangle \otimes \langle K'_\idot,B' \rangle$ of
a mixed complex $\langle \K_\idot,B \rangle$ in $\C$ and a mixed
complex $\langle K'_\idot,B' \rangle$ in $\C'$.

\begin{exa}\label{hh.mix}
For any ring $R$ and any object $E \in \Fun(\Lambda,R)$, the complex
$CH_\idot(E)$ with the differential $B$ of \eqref{B.ch} is an object
in $\D_{mix}(R)$.
\end{exa}

\begin{defn}
Assume that an abelian category $\C$ has countable products. The
{\em truncation} $\tr(\langle \K_\idot,B \rangle)$ of a mixed
complex $\langle \K_\idot,B \rangle$ in $\C$ is the product-total
complex of the bicomplex
$$
\begin{CD}
\K_\idot @>{B}>> \K_\idot[-1] @>{B}>> \K_\idot[-2] @>{B}>> ...
\end{CD}
$$
in $\C$.
\end{defn}

\proof[Proof of Lemma~\ref{B.kun}.]
For any mixed complex $\langle \K_\idot,B \rangle$ of $R$-modules,
let $\lambda(\langle \K_\idot,B \rangle)$ be the complex in
$\Fun(\Lambda,R)$ given by
$$
\lambda(\langle \K_\idot,B \rangle) = \tr(\langle \K_\idot,B \rangle
\otimes \K_\idot(R)),
$$
where $R$ in the right-hand side is shorthand for the constant
functor with value $R$. Then since $\K_\idot(R)$ is a complex of flat
objects in $\Fun(\Lambda,R)$, $\lambda$ respects quasiisomorphisms
and descends to a functor
\begin{equation}\label{mix.lambda}
\lambda:\D_{mix}(R) \to \D(\Lambda,R).
\end{equation}
Moreover, one checks easily that there exists a
quasiisomorphism
$$
\K_\idot(R) \cong \tr(\K_\idot(R) \otimes \K_\idot(R)),
$$
and this turns $\lambda$ into a triangulated tensor functor. The
Hochshild homology functor $CH_\idot(-)$ of Example~\ref{hh.mix} is
left-adjoint to $\lambda$. Since $\lambda$ is tensor, $CH_\idot(-)$
is pseudotensor by adjunction, so that we have a natural functorial
map
$$
CH_\idot(E \lotimes E') \to CH_\idot(E) \lotimes CH_\idot(E).
$$
This map is automatically compatible with the Connes-Tsygan
differentials, and it is inverse to the K\"unneth quasiisomorphism
\eqref{kunn}.
\endproof

\begin{remark}
In fact, $\lambda$ is a full embedding of tensor triangulated
categories, but we will not need this.
\end{remark}

\section{Cyclic homology of algebras.}\label{alg.sec}

Assume given a commutative ring $k$. To any associative unital
algebra $A$ over $k$ one associates a canonical object $A_{\hash}
\in \Fun(\Lambda,k)$ as follows:
\begin{itemize}
\item on objects, $A_{\hash}([n]) = A^{\otimes_k V([n])} = A^{\otimes_k
  n}$, with copies of $A$ numbered by elements $v \in V([n])$,
\item for any map $f:[n] \to [m]$, the map
$$
A_{\hash}(f):A^{\otimes_k n} = \bigotimes_{v \in V([m])}A^{\otimes_k f^{-1}(v)}
\to A^{\otimes m}
$$
is given by
\begin{equation}\label{a.hash.maps}
A_{\hash}(f) = \bigotimes_{v \in V([m])}m_{f^{-1}(v)},
\end{equation}
where $m_{f^{-1}(v)}:A^{\otimes_k f^{-1}(v)} \to A$ is the map which
multiplies the entries in the natural clockwise order on the vertices
$v' \in f^{-1}(v) \subset S^1$.
\end{itemize}
Assume that $A$ is flat as a $k$-module. Then by definition, the
{\em cyclic homology} $HC_\idot(A)$ is given by
$$
HC_\idot(A) = HC_\idot(A_{\hash}).
$$
The {\em Hochschild homology} $HH_\idot(A,M)$ with coefficients in
an $A$-bimodule $M$ is given by
$$
HH_\idot(A) = \Tor_\idot^{A^{opp} \otimes_k A}(A,M);
$$
explicitly, in degree $0$ we have
\begin{equation}\label{hh0}
HH_0(A,M) \cong M/\left\{am-ma\mid a \in A, m \in M \right\}.
\end{equation}
To simplify notation, one denotes $HH_\idot(A)=HH_\idot(A,A)$. One
easily shows using the bar resolution that we have a natural
identification
$$
HH_\idot(A) \cong HH_\idot(A_{\hash}).
$$
Under these identifications, \eqref{sp.se} reads as
\begin{equation}\label{hh.hc.A}
HH_\idot(A)[u] \Rightarrow HC_\idot(A).
\end{equation}
Note that by Lemma~\ref{la.ho}~\thetag{i}, one can equally well
replace $\Lambda$ with $\Lambda_n$ and $A_{\hash}$ with
$i_n^*A_{\hash}$; the resulting spectral sequence is the same.

In practice, the Hochschild homology groups $HH_\idot(A)$ can be
computed by the Hochschild homology complex $CH_\idot(A_{\hash})$,
denoted $CH_\idot(A)$ to simplify notation; its terms are given by
$$
CH_i(A) = A^{\otimes_k i+1},
$$
with a certain differential usually denoted by $b$. The Connes-Tsygan
differential $B$ can be lifted to the level of complexes, as in
\eqref{B.ch}. Explicitly, we have
\begin{equation}\label{B}
\begin{aligned}
B(a_0 \otimes \dots \otimes a_i) &= \sum_{j=0}^i 1 \otimes a_j \otimes
a_{j+1} \otimes \dots \otimes a_i \otimes a_0 \otimes \dots \otimes
a_{j-1}\\
&\quad - \sum_{j=0}^i a_j \otimes a_{j+1} \otimes \dots \otimes a_i
\otimes a_0 \otimes \dots \otimes a_{j-1} \otimes 1.
\end{aligned}
\end{equation}
Alternatively, one can fix an integer $n \geq 1$ and apply the
edgewise subdivision isomorphism $h_n$ of
Lemma~\ref{la.ho}~\thetag{i}. This gives a different complex
$CH_\idot^{(n)}(A) = CH_\idot(i_n^*A_{\hash})$ with terms
$$
CH^{(n)}_i(A) = A^{\otimes_k n(i+1)}
$$
computing the same Hochschild homology groups $HH_\idot(A)$. The
isomorphism $h_n$ can also be lifted to the level of complexes;
on $CH^{(n)}_i(A)$, it is given by
\begin{equation}\label{mu.l}
h_n(a_0^1 \otimes \dots \otimes a_i^1 \otimes \dots \otimes a_0^n
\otimes \dots \otimes a_i^n) =
a_0^1 \cdot \dots \cdot a_i^1 \cdot \dots
\cdot a_0^n \otimes a_1^n \otimes \dots \otimes a_i^n,
\end{equation}
or in words, ``multiply the first $(i+1)(n-1)+1$ terms, and keep the
last $i$ terms intact''.

In general, both $HH_\idot(A)$ and $HC_\idot(A)$ are just
$k$-modules. However, if $A$ is commutative, then $HH_\idot(A)$ becomes
commutative graded $k$-algebras, with the product induced by the
K\"unneth quasiisomorphism and the natural algebra map $A \otimes A
\to A$. By Lemma~\ref{B.kun}, the Connes-Tsygan differential $B$ is
a derivation with the respect to this algebra structure.

The main comparison theorem for Hochschild homology in the
commutative case is the classic result of Hochschild, Kostant, and
Rosenberg.

\begin{theorem}\label{hkr}
  Assume that $k$ is a field, the algebra $A$ is commutative and
  finitely generated over $k$, and that $X = \Spec A$ is
  smooth. Then we have natural isomorphisms
\begin{equation}\label{hkr.iso}
HH_i(A) \cong H^0(X,\Omega_X^i)
\end{equation}
for any $i \geq 0$, and these isomorphisms are compatible with
multiplication.\endproof
\end{theorem}

We want to emphasize that the Hochschild-Kostant-Rosenberg Theorem
requires no assumptions on $\cchar k$. As for cyclic homology, the
main result is the following.

\begin{theorem}\label{B.d.thm}
In the assumptions and under the identifications of
Theorem~\ref{hkr}, the Connes-Tsygan differential
$$
B:HH_i(A) \to HH_{i+1}(A),
$$
$i \geq 0$, becomes the de Rham differential $d:\Omega^i_X \to
\Omega^{i+1}_X$.
\end{theorem}

\proof{} This is extremely well-known if $\cchar k = 0$ (in fact, it
was this result which gave rise to cyclic homology as a separate
subject). For the convenience of the reader, let us show that the
result also holds when $\cchar k = p > 0$. Indeed, since $B$ is a
derivation, it suffices to check it on $HH_0(A) = A$ and $HH_1(A) =
\Omega^1(A)$, the module of K\"ahler differentials. But if $\cchar k
= p$, then for any $i$, $0 \leq i < p$ the
Hochschild-Kostant-Rosenberg isomorphism \eqref{hkr.iso} can be
lifted to $CH_i(A)$ and expressed by an explicit
formula
$$
P(a_0 \otimes \dots \otimes a_i) = \frac{1}{i!}a_0da_1 \wedge \dots
\wedge a_i.
$$
In particular, this is always possible in degrees $0$ and
$1$. Substituting this into \eqref{B} immediately proves the claim.
\endproof

\section{Trace maps.}\label{tra.sec}

The main technical tool in our approach to the Cartier map is the
projection \eqref{pi.n}; in this section, we analyze its homological
properties.

Fix a field $k$ and an integer $p \geq 1$. For a $k$-vector space
$M$ equipped with a representation of the cyclic group $\Z/p\Z$, let
$\sigma \in \Z/p\Z$ be the generator, and consider the spaces
$M^\sigma$, $M_\sigma$ of invariants and coinvariants. We have two
natural functorial maps between them: the map
\begin{equation}\label{e.l}
\e_p:M^\sigma \to M_\sigma
\end{equation}
obtained as the composition of the embedding $M^\sigma \subset M$
and the projection $M \to M_\sigma$, and the trace map
\begin{equation}\label{tr.l}
\tr_p = 1 + \sigma + \dots + \sigma^{p-1}:M_\sigma \to M^\sigma
\end{equation}
obtained by averaging over the group. We have
\begin{equation}\label{e.tr.l}
\tr_p \circ \e_p = p\id = \e_p \circ \tr_p.
\end{equation}
Moreover, say that a $k[\Z/p\Z]$-module is {\em induced} if $M = M'
\otimes_k k[\Z/p\Z]$ for some $k$-vector space $M'$. Then for an
induced $M$, the map $\tr_p$ is an isomorphism.

Now consider the projection $\pi^p:\Lambda_p \to \Lambda$ of
\eqref{pi.n}. As noted in Section~\ref{pre.sec}, $\pi^p$ is a
bifibration with fiber $\ppt_p = \ppt/(\Z/p\Z)$, so that for any $E
\in \Fun(\Lambda_p,R)$ and any $[n]$, $E([n])$ is a representation
of $\Z/p\Z$. Moreover, by \cite[Lemma 1.7]{ka0} we have
\begin{equation}\label{bc.pi}
\pi^p_!(E)([n]) \cong E([n])_\sigma, \qquad\pi^p_*(E)([n]) \cong
E([n])^\sigma.
\end{equation}
The maps $\e_p$ and $\tr_p$ are functorial in $[n]$ and give
two natural maps
\begin{equation}\label{e.tr.pi}
\e_p:\pi^p_*E \to \pi^p_!E, \qquad \tr_p:\pi^p_!E \to \pi^p_*E.
\end{equation}
If the $k[\Z/p\Z]$-module $E([n])$ is induced for every $[n] \in
\Lambda_p$, then the map $\tr_p$ of \eqref{e.tr.pi} is an
isomorphism, and $E$ is acyclic both for the right-exact functor
$\pi^p_!$ and for the left-exact functor $\pi^p_*$ (that is,
$L^i\pi^p_!E = \pi^p_!E$, $R^i\pi^p_*E = \pi^p_*E$). In particular,
for any $[n] \in \Lambda_p$, $\K_0([n]) \cong \K_1([n]) \cong
\Z[\Z/np\Z]$ are free modules over $\Z[\Z/p\Z]$, so that this
applies to functors of the form $\K_0(E) = E \otimes \K_0$, $\K_1(E)
= E \otimes \K_1$. Therefore we have a natural isomorphism
\begin{equation}\label{tr.k}
\tr_p:\pi^p_!\K_\idot(E) \cong \pi^p_*\K_\idot(E)
\end{equation}
for any $E \in \Fun(\Lambda_p,E)$. Moreover, since
$\pi^p_!\K_\idot(E) \cong L^\hdot\pi^p_!\K_\idot(E)$, we have
natural identifications
\begin{equation}\label{pi.k}
HC_\idot(\K_\idot(E)) \cong HC_\idot(\pi^p_!\K_\idot(E)) \cong
HC_\idot(\pi^p_*\K_\idot(E)).
\end{equation}
By Lemma~\ref{la.ho}~\thetag{ii}, all these groups are further
identified with $HH_\idot(E)$. The homology of the complex
$\pi^p_!\K_\idot(E)$ in degree $0$ coincides with $\pi^p_!E$, with the
identification provided by the map $\pi^p_!(\kappa_1):\pi^p_!\K_0(E)
\to \pi^p_!E$, while by \eqref{tr.k}, the homology of the same
complex in degree $1$ is $\pi^p_*E$, with the identification
provided by the map $\pi^p_*(\kappa_1):\pi^p_*E \to \pi^p_*\K_1(E)
\cong \pi^p_!\K_1(E)$. Moreover, denote by
\begin{equation}\label{B.tr}
\wt{B}:\pi^p_!\K_0(E) \to \pi^p_!\K_1(E)
\end{equation}
the map obtained as the composition
$$
\begin{CD}
\pi^p_!\K_0(E) @>{\pi^p_!(\kappa_0)}>> \pi^p_!E @>{\tr_p}>> \pi^p_*E
@>{\pi^p_*(\kappa_1)}>> \pi^p_*\K_1(E) @>{\tr_p^{-1}}>> \pi^p_!\K_1(E),
\end{CD}
$$
where $\tr_p^{-1}$ is the inverse to the isomorphism
\eqref{tr.k}. Then by the functoriality of the trace map $\tr_p$, we
have $\tr_p \circ \pi^p_!(\kappa_1) = \pi^p_*(\kappa_1) \circ
\tr_p$, so that
\begin{equation}\label{B.B}
\wt{B} = \tr_p^{-1} \circ \pi^p_*(\kappa_1) \circ \tr_p \circ
\pi^p_!(\kappa_0) = \pi^p_!(\kappa_1) \circ \pi^p_!(\kappa_0) =
\pi^p_!(B),
\end{equation}
where $B:\K_1(E) \to \K_0(E)$ is the natural map \eqref{B.0}. Under
the identifications \eqref{pi.k}, the map $\wt{B}$ of \eqref{B.tr}
then induces the Connes-Tsygan differential \eqref{c.ts} on
$HH_\idot(E)$.

\begin{lemma}\label{k.l.l}
For any $p \geq 1$ and any $E \in \Fun(\Lambda,k)$, we have
natural isomorphisms of complexes
$$
\K_\idot(E) \cong \pi^p_!\K_\idot(\pi^{p*}E) \cong
\pi^p_*\K_\idot(\pi^{p*}E).
$$
\end{lemma}

\proof{} By the projection formula \cite[Lemma 1.7]{ka0}, we have
$\pi^p_!\K_\idot(\pi^{p*}E) \cong E \otimes \pi^p_!\pi^{p*}\Z$ and
$\pi^p_*\K_\idot(\pi^{p*}E) \cong E \otimes \pi^p_*\pi^{p*}\Z$, so
that it suffices to consider the case $E = \Z$. In this case, the
claim immediately follows from the definition of the objects $\K_0$
and $\K_1$.
\endproof

Next, consider the derived functor $L^\hdot\pi^p_!$. We first
observe the following. Since $\pi^p$ is a bifibration, \cite[Lemma
  1.7]{ka0} shows that for any $E \in \Fun(\Lambda_p,k)$ and any $n
\geq 1$, we also have a natural identification
$$
L^\hdot\pi^p_!(E)([n]) \cong H_\idot(\Z/p\Z,E([n])),
$$
a derived version of \eqref{bc.pi}, where the group $\Z/p\Z$
acting on $E([n])$ is generated by the automorphism $\sigma$. In
particular, this applies to the constant functor $k \in
\Fun(\Lambda_p,k)$. But the homology of any group with constant
coefficients is a coalgebra, and the homology of the group with any
coefficients is a comodule over this coalgebra. We claim that for
the group $\Z/p\Z$, this can be made to work relatively over the
category $\Lambda$.

To see this, let $F_\idot \in \Fun(\Lambda_p,k)$ be any resolution
of the constant functor $k$ by functors $F_i$ such that $F_i([n])$
is a free $k[\Z/p\Z]$-module for any $i \geq 0$, $n \geq 1$ -- for
example, we can use the resolution \eqref{reso}. Then for any $E \in
\Fun(\Lambda_p,k)$ and any $i \geq 0$, the product $F_i \otimes E$
has the same property, so that the resolution $F_\idot \otimes E$
can be used to compute $L^\hdot\pi_{p!}(E)$. Denote
$$
\wt{\Lambda}_p = \Lambda_p \times_\Lambda \Lambda_p,
$$
let $\delta:\Lambda_p \to \wt{\Lambda}_p$ be the diagonal embedding,
and let $\wt{\pi}^p:\wt{\Lambda}_p \to \Lambda$ be the natural
projection. Then for any $E \in \Fun(\Lambda_p,k)$, we have
$$
\begin{aligned}
L^\hdot\wt{\pi}^p_!(k \boxtimes E) &\cong \wt{\pi}^p_!(F_\idot
\boxtimes (F_\idot \otimes E)) \cong \\
&\cong \pi^p_!(F_\idot) \otimes
\pi^p_!(F_\idot \otimes E) \cong L^\hdot\pi^p_!(k) \otimes
L^\hdot\pi^p_!(E).
\end{aligned}
$$
On the other hand $\delta^*(k \boxtimes E) \cong E$, and the natural
projection
$$
\pi^p_!\delta^*(F_\idot \boxtimes (F_\idot \otimes E)) =
\pi^p_!(F_\idot \otimes F_\idot \otimes E) \to
\wt{\pi}^p_!(F_\idot \boxtimes (F_\idot \otimes E))
$$
induces a map
\begin{equation}\label{wt.a}
\wt{a}:L^\hdot\pi^p_!(E) \to L^\hdot\pi^p_!(k) \otimes
L^\hdot\pi^p_!(E).
\end{equation}
This is our coaction map.

\medskip

Now assume that $p$ is a prime and $\cchar k = p$. Then we have the
following splitting result.

\begin{lemma}\label{p.spl.le}
For any $E \in \Fun(\Lambda,k)$, there exists a natural isomorphism
\begin{equation}\label{spl.pi}
L^\hdot\pi^p_!\pi^{p*}E \cong \K_\idot(E)[u] = \bigoplus_{n \geq
  0}\K_\idot(E)[2n].
\end{equation}
\end{lemma}

\proof{} To compute $L^\hdot\pi^p_!(E)$, we can use the resolution
\eqref{reso}. Then by \eqref{B.B}, the differential $\pi^p_!(B)$ in
the resulting complex is given by the trace map
$$
\tr_p = 1 + \sigma + \dots + \sigma^{p-1}:\pi^p_!\pi^{p*}E \to
    \pi^p_*\pi^{p*}E
$$
of \eqref{e.tr.pi}. But $\sigma = \id$ on $\pi^{p*}E$, so that the
map $\tr_p$ is just multiplication by $p$, thus equal to $0$ by
assumption. Combining this with the canonical isomorphism of
Lemma~\ref{k.l.l}, we obtain the desired isomorphism already on the
level of complexes.
\endproof

Composing the coaction map $\wt{a}$ of \eqref{wt.a} with the
projection $L^\hdot\pi^p_!(k) \to \K_\idot(k)$ onto the first
summand in \eqref{spl.pi}, we obtain a canonical functorial map
\begin{equation}\label{nowt.a}
a:L^\hdot\pi^p_!(E) \to L^\hdot\pi^p_!(E) \otimes \K_\idot(k)
\end{equation}
in the derived category $\D(\Lambda,k)$. 

To proceed further, we need to assume that $p$ is odd, and impose a
condition on the object $E \in \Fun(\Lambda_p,k)$. Namely, recall
that if $p \neq 2$, the cohomology algebra $H^\hdot(\Z/p\Z,k)$ is
given by
$$
H^\hdot(\Z/p\Z,k) \cong k[u]\langle \eps \rangle,
$$
where in the right-hand side, we have the free graded-commutative
algebra on one generator $u$ of degree $2$ and one generator $\eps$
of degree $1$. For any $k[\Z/p\Z]$-module $E$, multiplication by $u$
induces an isomorphism
$$
H_{i+2}(\Z/p\Z,E) \cong H_i(\Z/p\Z,E)
$$
for any $i \geq 1$.

\begin{defn}\label{tight}
A $k[\Z/p\Z]$-module $E$ is {\em tight} if the map
\begin{equation}\label{eps.eq.0}
H_{2i+1}(\Z/p\Z,E) \to H_{2i}(\Z/p\Z,E)
\end{equation}
given by multiplication by $\eps \in H^1(\Z/p\Z,k)$ is an
isomorphism for any $i \geq 1$. An object $E \in \Fun(\Lambda_p,k)$
is {\em tight} if $E([n])$ is a tight $k[\Z/p\Z]$-module for any
$[n] \in \Lambda_p$.
\end{defn}

To see the tightness condition explicitly, note that by
\eqref{e.tr.l} we have $\e_p \circ \tr_p = 0$ for every
$k[\Z/p\Z]$-module $M$, so that $\e_p$ induces a map
\begin{equation}\label{eps.eq}
\Coker \tr_p \to \Ker \tr_p.
\end{equation}
If one computes the homology groups $H_\idot(\Z/p\Z,E)$ by the
standard periodic complex, then $\Ker \tr$ resp.\ $\Coker \tr$ is
identified with $H_{2\idot + 1}(\Z/p\Z,E)$
resp.\ $H_{2\idot}(\Z/p\Z,E)$, and the map \eqref{eps.eq} is exactly
the map \eqref{eps.eq.0}.

\begin{exa}\label{tight.exa}
A trivial $k[\Z/p\Z]$-module $k$ is tight; so is a free module
  $k[\Z/p\Z]$.
\end{exa}

\begin{remark}\label{Efi}
Example~\ref{tight.exa} is essentially exhaustive: all
indecomposable tight $k[\Z/p\Z]$-modules are of this form.
\end{remark}

We note that since $\eps^2=0$, for any tight $k[\Z/p\Z]$-module $E$,
the map
\begin{equation}\label{2i}
H_{2i}(\Z/p\Z,E) \to H_{2i-1}(\Z/p\Z,E)
\end{equation}
given by multiplication by $\eps$ is equal to $0$. We also note that
for a tight $E \in \Fun(\Lambda_p,k)$, all the homology objects of
the direct image $L^\hdot\pi_{p!}(E)$ in degree $\geq 1$ are
canonically identified; we will denote this homology object by
\begin{equation}\label{I.E}
\I(E) \in \Fun(\Lambda,k).
\end{equation}
We can now prove our main result in this section. Consider the
standard $t$-structure on the derived category $\D(\Lambda,k)$, and
let $\tau_{\leq \idot}$, $\tau_{[\idot,\idot]}$ be the corresponding
truncation functors, so that for any integers $m \geq n$ and any $E
\in \D(\Lambda,k)$, we have a natural exact triangle
$$
\begin{CD}
\tau_{\leq n}(E) @>>> \tau_{\leq m}(E) @>>> \tau_{[n+1,m]}(E) @>>>
\end{CD}
$$
in the category $\D(\Lambda,k)$.

\begin{lemma}\label{gr.le}
Assume given a tight object $E \in \Fun(\Lambda_p,E)$, and assume
that $p$ is odd. Then for any $i \geq 1$, the coaction map $a$ of
\eqref{nowt.a} induces a map
$$
\tau_{\leq -2i}L^\hdot\pi^p_!(E) \to \tau_{\leq
  -2i}L^\hdot\pi^p_!(E) \otimes \K_\idot(k),
$$
and the composition map
$$
\begin{CD}
\tau_{[-2i-1,-2i]}L^\hdot\pi^p_!(E) @>>>
\tau_{[-2i-1,-2i]}L^\hdot\pi^p_!(E) \otimes \K_\idot(k) @>>>\\
@>>> L^{2i}\pi^p_!(E) \otimes \K_\idot(k) \cong
\K_\idot(\I(E))
\end{CD}
$$
is a quasiisomorphism.
\end{lemma}

\proof{} To prove the first claim, it suffices to check that the
composition map
$$
\tau_{\leq -2i}L^\hdot\pi^p_!(E) \to L^\hdot\pi^p_!(E) \otimes
\K_\idot(k) \to \tau_{\geq -2i+1}L^\hdot\pi^p_!(E) \otimes
\K_\idot(k)
$$
is equal to $0$. The left-hand side is in $\D_{\leq -2i}(\Lambda,k)$
with respect to the standard $t$-structure, and the right-hand side
is in $\D^{\geq -2i}(\Lambda,k)$; therefore it suffices to check
that the map is equal to $0$ on homology objects of degree
$2i$. This can be checked after evaluating at every $[n] \in
\Lambda$, and the corresponding map is the map \eqref{2i} which is
indeed equal to $0$ for a tight $E$. Analogously, to prove the
second claim, it suffices to check that the map is an isomorphism in
homological degrees $2i$ and $2i+1$; this follows immediately from
the definition of a tight object.
\endproof

\section{Cartier isomorphism.}\label{car.sec}

We can now present the construction of our generalized Cartier
isomorphism. We start with some linear algebra. Assume that our
field $k$ of positive characteristic $p = \cchar k > 0$ is perfect,
and let $M$ be a $k$-vector space. Let the group $\Z/p\Z$ act on the
$p$-th tensor power $M^{\otimes p}$ so that the generator $\sigma
\in \Z/p\Z$ acts by the cyclic permutation of order $p$. Consider
the trace map
\begin{equation}\label{tr.M}
\tr_p:\left(M^{\otimes p}\right)_\sigma \to \left(M^{\otimes
  p}\right)^\sigma
\end{equation}
of \eqref{tr.l}.

\begin{lemma}\label{p-th.le}
\begin{enumerate}
\item The $k[\Z/p\Z]$-module $M^{\otimes p}$ is tight in the sense
  of Definition~\ref{tight}.
\item Sending $m \in M$ to $m^{\otimes p}$ gives a well-defined
  additive Frobenius-semi\-li\-near map
$$
\psi:M \to \left(M^{\otimes p}\right)_\sigma
$$
which identifies $M$ with the kernel of the map \eqref{tr.M}. The
dual map
$$
\wh{\psi}:\left(M^{\otimes p}\right)^\sigma \to M
$$
identifies $M$ with the cokernel of map $\tr$.
\item For any $n \geq 0$, we have a commutative diagram
$$
\begin{CD}
\left(M^{\otimes p^n}\right)^\sigma @>{\wh{\psi}^n}>> M\\
@VVV @VV{\psi^n}V\\
M^{\otimes p^n} @>>> \left(M^{\otimes p^n}\right)_\sigma,
\end{CD}
$$
where the left vertical map is the natural embedding, and the bottom
map is the natural projection.
\end{enumerate}
\end{lemma}

\proof{} To prove \thetag{i}, choose a basis in $M$, so that $M =
k[S]$ for some set $S$. Then
\begin{equation}\label{splo}
M^{\otimes p} = k[S^p] = k[S] \oplus k[S^p \setminus S],
\end{equation}
where $S$ is embedded into $S^p$ as the diagonal. The first summand
is trivial, and the second summand is free, and both are tight.

\thetag{ii} is actually \cite[Lemma 2.3]{ka0}, but let us reproduce
the proof for the convenience of the reader. Consider the standard
periodic complex
$$
\begin{CD}
@>{\id + \sigma + \dots + \sigma^{p-1}}>> M^{\otimes p} @>{\id -
  \sigma}>> M^{\otimes p} @>{\id + \sigma + \dots + \sigma^{p-1}}>>
\end{CD}
$$
computing the Tate homology $\vH_\idot(\Z/p\Z,M^{\otimes
  p})$. Define a non-additive map $\psi:M \to M^{\otimes p}$ by
$\psi(m) = m^{\otimes p}$. Then this map goes into the kernel of the
differential of the complex (both in odd and even degrees), and it
becomes additive after we project onto the first summand in
\eqref{splo}. Moreover, since $k$ is perfect, it becomes a
quasiisomorphism (again both in odd and in even degrees). It remains
to notice that the second summand in \eqref{splo} has no Tate
homology, so that the corresponding complex is acyclic, and $\psi$
is additive module the image of the differential. Therefore the
induced map
$$
\psi:M \to \Coker(\id - \sigma) = \left(M^{\otimes p}\right)_\sigma
$$
is additive and identifies $M$ with the kernel of $\tr$, as
required. The second claim follows by duality.

Finally, for \thetag{iii}, again choose a basis $S$ in $M$, so that
$S^{p^n}/(\Z/p^n\Z)$ gives a natural basis both in $\left(M^{\otimes
  p^n}\right)^\sigma$ and in $\left(M^{\otimes
  p^n}\right)_\sigma$, and note that the composition of the natural
embedding and the natural projection is a diagonal operator in this
basis, with the entry correponding to a $\Z/p^n\Z$-orbit $S' \subset
S^{p^n}$ equal to its cardinality $|S'|$. Since $\cchar k = p$, this
is $1$ for the one-point orbits, and $0$ otherwise.
\endproof

Let now $A$ be an associative algebra over $k$, and consider the
corresponding object $A_{\hash} \in \Fun(\Lambda,k)$ of
Section~\ref{alg.sec} and its pullback $i_p^*(A_{\hash}) \in
\Fun(\Lambda_p,k)$ with respect to the embedding $i_p:\Lambda_p \to
\Lambda$. Note that the canonical maps $\psi$ and $\wh{\psi}$ of
Lemma~\ref{p-th.le} fit together to give canonical maps
\begin{equation}\label{psi.la}
\psi:A_\hash \to \pi^p_!i_p^*A_\hash, \qquad
\wh{\psi}:\pi^p_*i_p^*A_\hash \to A_\hash.
\end{equation}

\begin{corr}\label{I.A}
The object $i_p^*(A_\hash) \in \Fun(\Lambda_p,k)$ is tight in the
sense of Definition~\ref{tight}, and the map $\psi$ of
\eqref{psi.la} induces a Frobenius-semilinear isomorphism
$$
A_\hash \cong \I(i_p^*(A_\hash)).
$$
\end{corr}

\proof{} Immediately follows from Lemma~\ref{p-th.le}. \endproof

To simplify notation, for any $E \in \Fun(\Lambda,k)$ and any $n
\geq 1$, let us denote by $\K^n_\idot(E)$ the length-$2$ complex in
$\Fun(\Lambda,k)$ given by
$$
\K^n_\idot(E) = \pi^n_!\K_\idot(i_n^*E) \cong
\pi^n_*\K_\idot(i_n^*E).
$$
Then \eqref{pi.k} together with the isomorphism $h_n$ of
Lemma~\ref{la.ho}~\thetag{i} provide a canonical isomorphism
\begin{equation}\label{K.l.cor}
HH_\idot(E) \cong HC_\idot(\K^n_\idot(E))
\end{equation}
for any $E \in \Fun(\Lambda,k)$ and any $n \geq 1$. In particular,
for $E = A_{\hash}$ and $n=p$, we have a natural identification
$$
HC_\idot(\K^p_\idot(A_{\hash})) \cong HH_\idot(A_{\hash}) \cong
HH_\idot(A).
$$

\begin{defn}\label{B.Z.Kp}
Assume given an object $E \in \Fun(\Lambda,k)$. Then the
subcomplexes $B\K^p_\idot(E),Z\K^p_\idot(E) \subset
\K^p_\idot(E)$ are the image resp. the kernel of the map
$$
\wt{B}:\K^p_\idot(E) \to \K^p_\idot(E)[-1]
$$
of \eqref{B.tr}. For any associative $k$-algebra $A$, we denote
$$
\begin{aligned}
BHH_\idot(A) &= HC_\idot(B\K^p_\idot(A_{\hash})),\\
ZHH_\idot(A) &= HC_\idot(Z\K^p_\idot(A_{\hash})).
\end{aligned}
$$
\end{defn}

Alternatively, by Lemma~\ref{p-th.le}~\thetag{ii} and \eqref{B.B}, we
could define $B\K^p_\idot(E)$ and $Z\K^p_\idot(E)$ by saying that on
the level of homology, we have
$$
H_i(B\K^p_\idot(E))=
\begin{cases}
0, &\quad i = 0,\\
\Ker \wh{\psi} \subset \pi^{p,1}_*i_p^*A_\hash = H_1(\K^p_\idot(E)),
&\quad i=1,
\end{cases}
$$
and
$$
H_i(Z\K^p_\idot(E)) = 
\begin{cases}
\Im \psi \subset \pi^{p,1}_!i_p^*A_\hash = H_0(\K^p_\idot(E)),
&\quad i=0,\\
\pi^{p,1}_*i_p^*A_\hash = H_1(\K^p_\idot(E)), &\quad i=1,
\end{cases}
$$
where $\psi$ and $\wh{\psi}$ are the canonical maps
\eqref{psi.la}. In any case, we have
$$
B\K^p_\idot(E) \subset Z\K^p_\idot(E) \subset
\K^p_\idot(E),
$$
and these embeddings induce natural maps
\begin{equation}\label{B.Z}
\begin{CD}
BHH_\idot(A) @>{\xi}>> ZHH_\idot(A) @>{\zeta}>> HH_\idot(A).
\end{CD}
\end{equation}
In the general case, I do not know whether these maps are injective
or not. One obvious observation is that the Connes-Tsygan
differential \eqref{c.ts}, being induced by the map $\wt{B}$ of
\eqref{B.tr}, factors through a map
\begin{equation}\label{b.z}
\beta:HH_{\idot} \to BHH_{\idot+1}(A),
\end{equation}
and this map actually fits into a long exact sequence
\begin{equation}\label{z.lg}
\begin{CD}
ZHH_\idot(A) @>{\zeta}>> HH_\idot(A) @>{\beta}>> BHH_{\idot+1}(A) @>>>
\end{CD}
\end{equation}
This has the following corollary. Denote by
$zHH_\idot(Z),bHH_\idot(A) \subset HH_\idot(A)$ the kernel resp. the
image of the Connes-Tsygan differential, and denote by
$z'HH_\idot(Z),b'HH_\idot(A) \subset HH_\idot(A)$ the image of the
map $\zeta$ resp. $\zeta \circ \xi$. Then we have inclusions
\begin{equation}\label{incl}
bHH_\idot(A) \subset b'HH_\idot(A) \subset z'HH_\idot(A) \subset
zHH_\idot(A).
\end{equation}

\begin{lemma}\label{car.le}
Assume that $p$ is odd. Then there exists a natural
Frobenius-semilinear identification
$$
HC_\idot(Z\K^p_\idot(A_{\hash})/B\K^p_\idot(A_{\hash})) \cong
HH_\idot(A).
$$
\end{lemma}

\proof{} Let us compute the direct image
$L^\hdot\pi^p_!i_p^*A_{\hash}$ by the standard resolution
\eqref{reso}. Then we obtain a natural quasiisomorphism
$$
Z\K^p_\idot(A_{\hash})/B\K^p_\idot(A_{\hash})) \cong 
\tau_{[-2i-1,-2i]}L^\hdot\pi^p_!i_l^*A_{\hash}.
$$
Combining this with Lemma~\ref{gr.le} and Corollary~\ref{I.A}, we
obtain the desired identification.
\endproof

\begin{corr}
For any associative algebra $A$ over a perfect field $k$ of
characteristic $p \neq 2$, there exists a canonical long exact
sequence
\begin{equation}\label{car.seq}
\begin{CD}
BHH_\idot(A) @>{\xi}>> ZHH_\idot(A) @>{C}>> HH_\idot(A) @>>>,
\end{CD}
\end{equation}
where $\xi$ one of the canonical maps \eqref{B.Z}, and the map $C$
is induced by the isomorphism of Lemma~\ref{car.le}.
\end{corr}

\proof{} Clear. \endproof

\begin{defn}\label{car.def}
The {\em non-commutative Cartier map} for the algebra $A$ is the
canonical map $C$ of \eqref{car.seq}.
\end{defn}

\section{Comparison.}\label{comp.sec}

We now need to justify the name ``Cartier map'' used in
Definition~\ref{car.def}. The comparison result is the following.

\begin{prop}\label{car.comp}
Assume that the algebra $A$ is commutative and finitely generated
over the perfect field $k$, $\cchar k \neq 2$, and that its spectrum
$X = \Spec A$ is smooth. Then the canonical maps $\xi$, $\zeta$ of
\eqref{B.Z} are injective, the map $\beta$ of \eqref{b.z} is
surjective, the connecting differential in the long exact sequence
\eqref{car.seq} is trivial, the Hochschild-Kostant-Rosenberg
isomorphism induces isomorphisms
$$
BHH_i(A) \cong B\Omega^i(X), \qquad ZHH_i(A) \cong Z\Omega^i(X)
$$
for every $i \geq 0$, and the non-commutative Cartier map of
Definition~\ref{car.def} coincides with the classical Cartier
isomorphism.
\end{prop}

For the proof, we first assume given a monoid $G$, and consider its
group algebra $A = k[G]$. Then the diagonal map $G \to G^p$ induces
a $\Z/p\Z$-invariant algebra map
$$
\phi:A \to A^{\otimes p},
$$
and this map together with its tensor powers defines a map
$\phi:\pi^{p*}A_\hash \to i_p^*A_\hash$
and a corresponding map 
$$
\phi:\pi^p_!\K_\idot(\pi^{p*}A_\hash) \to
\pi^p_!\K_\idot(i_p^*A_\hash).
$$
By Lemma~\ref{p.spl.le}, at the level of the derived category, this
can be rewritten as
\begin{equation}\label{Phi.eq}
\Phi:\K_\idot(A_\hash) \to \K_\idot^p(A_\hash)
\end{equation}
(in fact, for any $E \in \Fun(\Lambda,k)$, we have
$\pi^p_!\K_\idot(\pi^{p*}E) \cong \K_\idot(E)$ on the nose, but we
will not need this).

\begin{lemma}\label{grp.le}
The map $\Phi$ of \eqref{Phi.eq} factors through a map
$$
\Phi:\K_\idot(A_\hash) \to Z\K_\idot^p(A_\hash) \subset
\K_\idot^p(A_\hash),
$$
and the induced composition map
$$
\begin{CD}
HH_\idot(A) @>{\Phi}>> ZHH_\idot(A) @>{C}>> HH_\idot(A)
\end{CD}
$$
is equal to the identity map.
\end{lemma}

\proof{} For any vector space $M = k[S]$ with a base $S$, the
diagonal map $S \to S^p$ induces a map $\phi:M \to M^{\otimes p} =
k[S^p]$. This map obviously lands in the subspace $(M^{\otimes
  p})^\sigma \subset M^{\otimes p}$, and its composition with the
canonical map $\wh{\psi}$ of \eqref{psi.la} is equal to the
identity. Applying this pointwise to the object $A_\hash$, we see
that the composition
$$
\begin{CD}
A_\hash @>{\phi}>> \pi^p_*i_p^*A_\hash @>{\wh{\psi}}>>
\I(i_p^*A_\hash) \cong A_\hash
\end{CD}
$$
is also equal to the identity. Applying the identification of
Lemma~\ref{gr.le}, we obtain the claim.
\endproof

This immediately implies that in the case $A = k[G]$, the connecting
differential in the long exact sequence \eqref{car.seq} is equal to
$0$. We in fact have
\begin{equation}\label{xi.phi}
ZHH_\idot(A) \cong \xi(BHH_\idot(A)) \oplus \Phi(HH_\idot(A)),
\end{equation}
and both $\xi$ and $\Phi$ are injective.

\medskip

Next, fix an integer $n \geq 1$, and consider the algebra 
$$
A_n = k[t_1,t_1^{-1},\dots,t_n,t_n^{-1}]
$$
of Laurent polynomials in $n$ variables. Then $X = \Spec A_n \cong
\mathbb{G}_m^n$ is smooth, and $A_n$ satisfies the assumptions of
the Hochschild-Kostant-Rosenberg Theorem and
Theorem~\ref{B.d.thm}. The algebra $\Omega^\hdot(X)$ of differential
forms on $X$ is the free graded-commutative $A_n$-algebra generated
by logarithmic derivatives $dt_i/t_i$, $1 \leq i \leq n$. The
easiest way to see the classical Cartier map for $X$ is to consider
the multiplicative Frobenius-semilinear map
\begin{equation}\label{C.-1}
C^{-1}:\Omega^\hdot(X) \to \Omega^\hdot(X)
\end{equation}
given by $C^{-1}(t_i) = t_i^p$, $C^{-1}(dt_i/t_i)=dt_i/t_i$, $1 \leq
i \leq n$. Then $C^{-1}$ actually maps $\Omega^\hdot(X)$ into
$Z^\hdot(X) \subset \Omega^\hdot(X)$, and its composition with the
classical Cartier map $Z^\hdot(X) \to \Omega^\hdot(X)$ is equal to
the identity.

On the other hand, we have $A_n = k[\Z^n]$, the group algebra of the
free abelian group on $n$ generators. Therefore Lemma~\ref{grp.le}
also applies to the algebra $A_n$, and in particular, we have a
natural map
$$
\Phi:HH_\idot(A_n) \to ZHH_\idot(A_n).
$$
Moreover, fix another integer $l$, $0 \leq l \leq n$, and consider
the subalgebra 
\begin{equation}\label{a.ln}
A_{l,n} = k [t_1,\dots,t_l,t_{l+1},t_{l+1}^{-1},\dots,t_n,t_n^{-1}]
  \subset A_n
\end{equation}
of Laurent polynomials having no poles in the first $l$
variables. Then $Y = \Spec A_{m,n}$ is isomorphic to $\mathbb{G}_a^l
\times \mathbb{G}_m^{n-l}$, thus also satisfies the assumptions of
the comparison theorems, and the map $C^{-1}$ sends $\Omega^\hdot(Y)
\subset \Omega^\hdot(X)$ into itself. On the other hand, $A_{l,n}
\cong k[\N^l \times \Z^{n-l}]$, where $\N$ is the monoid of
non-negative integers, so that the map $\Phi$ is well-defined.

\begin{lemma}\label{phi.car}
Under the Hochschild-Kostant-Rosenberg isomorphism \eqref{hkr.iso},
the composition $\zeta \circ \Phi:HH_\idot(A_{l,n}) \to
HH_\idot(A_{l,n})$ is identified with the map $C^{-1}$ of
\eqref{C.-1}.
\end{lemma}

\proof{} Both sides are compatible with the K\"unneth isomorphism,
so that it suffices to consider the case $n=1$, $l= 0,1$. Moreover,
the natural maps $HH_\idot(A_{1,1}) \to HH_\idot(A_{0,1})$ are
injective, and maps $C^{-1}$ and $\Phi$ for the algebra $A_{1,1}$
are obtained by restriction from the corresponding maps for the
algebra $A_1=A_{0,1}$. Thus is suffices to consider the algebra $A =
A_1 = k[t,t^{-1}] = k[\Z]$. The Hochschild homology $HH_\idot(A)$ is
then only non-trivial in degrees $0$ and $1$. The $i$-th term
$CH_i(A)$ of the Hochschild complex of the algebra $A$ is given by
$$
CH_i(A) = A^{\otimes i+1} = k[\Z^{i+1}],
$$
and using the explicit form \eqref{mu.l} of the isomorphism $h_l$
of Lemma~\ref{la.ho}, one immediately checks that on $CH_i(A)$, the
composition $\zeta \circ \Phi$ is given by
$$
\zeta(\Phi(\langle g_0,\dots,g_i\rangle)) = \langle
g_0 + (p-1)(g_0 + \dots + g_i),g_1,\dots,g_i\rangle
$$
for any basis element $\langle g_0,\dots,g_i\rangle \in
\Z^{i+1}$. For $i=0$ and $i=1$, this reads as
$$
\begin{aligned}
\zeta(\Phi(t^n)) &= t^{pn},\\
\zeta(\Phi(\langle t^{n_0},t^{n_1}\rangle)) &= \langle
t^{pn_0+(p-1)n_1}, t^{n_1}\rangle.
\end{aligned}
$$
This obviously coincides with $C^{-1}$ in degree $0$, and since the
Hochschild-Kostant-Rosenberg isomorphism sends $\langle
t^{n_0},t^{n_1} \rangle$ to $n_1t^{n_1+n_0-1}dt$, we also get the
result in degree $1$.
\endproof

\begin{lemma}\label{car.tor}
For any $n \geq 1$, $0 \leq l \leq n$, Proposition~\ref{car.comp}
holds for the algebra $A=A_{l,n}$ of \eqref{a.ln}.
\end{lemma}

\proof{} The classical Cartier isomorphism and the inverse map
$C^{-1}$ give a splitting
$$
zHH_\idot(A) = bHH_\idot(A) \oplus C^{-1}(HH_\idot(A)),
$$
and by Lemma~\ref{phi.car}, $C^{-1}$ factors through $\zeta$, so
that by \eqref{incl}, we have $zHH_\idot(A) =
z'HH_\idot(A)$. Moreover, since $C^{-1}$ is injective, $\zeta$ is
injective on the direct summand $\Phi(HH_\idot(A))$ of the
decomposition \eqref{xi.phi}. Assume that the natural map
$\beta:HH_i(A) \to BHH_{i+1}(A)$ of \eqref{b.z} is surjective for $i
< l$ for some integer $l \geq 0$. Then we have
\begin{equation}\label{B.b}
BHH_i(A) \cong HH_{i-1}(A)/zHH_{i-1}(A) \cong bHH_i(A)
\end{equation}
for $i \leq l$, so that $\zeta$ is injective on $\xi(BHH_i(A))$ and
$bHH_i(A) = b'HH_i(A)$ in this range of degrees. Thus for $i \leq l$,
$\zeta$ is injective on the whole $ZHH_i(A)$, and the connecting
homomorphism in \eqref{z.lg} vanishes, so that $\beta$ is surjective
on $HH_l(A)$, too.

By induction, $\xi$ is injective and $\beta$ is surjective in all
degrees, and we have $BHH_\idot(A) \cong bHH_\idot(A)$,
$ZHH_\idot(A) \cong zHH_\idot(A)$. It remains to check that the
classical Cartier map and the non-commutative Cartier map give the
same identification $zHH_\idot(A)/bHH_\idot(A) \cong
HH_\idot(A)$. This immediately follows from Lemma~\ref{phi.car},
since they have the same inverse.
\endproof

\proof[Proof of Proposition~\ref{car.comp}.] By assumption, $A$ is
finitely generated, $X = \Spec A$ is smooth, and $HH_\idot(A) =
\Omega^\hdot(X)$ is spanned by differential forms $f_0 df_1 \wedge
\dots \wedge df_i$, $f_0,\dots,f_i \in A$. Every such form is the
image of the standard form
$$
\tau_i = t_0 dt_1 \wedge \dots \wedge dt_i
$$
under the map $\alpha:A_{i,i} = k[t_0,\dots,t_i] \to A$ sending
$t_j$ to $f_j$, $0 \leq j \leq i$. By Lemma~\ref{car.tor}, this
immediately implies that the non-commutative Cartier map
$C:ZHH_\idot(A) \to HH_\idot(A)$ is surjective, so that
$\xi:BHH_\idot(A) \to ZHH_\idot(A)$ is injective. Moreover, fix an
element $\wt{\tau}_i \in ZHH_i(A_{i,i})=zHH_i(A_{i,i}) \subset
HH_i(A_{i,i})$ such that $C(\wt{\tau}_i) = \tau_i$. Then by the
classic Cartier isomorphism, every closed form $a \in Z^i(X) \cong
zHH_i(A)$ can be expressed as a finite sum
$$
a = \overline{a} + \sum_j \alpha_j(\wt{\tau}_i),
$$
where $\overline{a} \in B^i(X) \cong bHH_i(A)$ is an exact form,
and $\alpha_j:A_{i,i} \to A$ are some algebra maps. Since $\wt{\tau}_i$
lies in $z'HH_i(A_{i,i})=zHH_i(A_{i,i})$, $a$ must also must lie in
$z'HH_i(A)$ by \eqref{incl}, so that in fact $z'HH_\idot(A) =
zHH_\idot(A)$.

Now, as in the proof of Lemma~\ref{car.tor}, assume by induction
that the natural map $\beta:HH_i(A) \to BHH_{i+1}(A)$ of \eqref{b.z}
is surjective for $i < l$ for some integer $l \geq 0$. Then again,
we have \eqref{B.b} for $i \leq l$, so that $\zeta$ is injective on
$\xi(BHH_i(A))$ and $bHH_i(A) = b'HH_i(A)$ in this range of
degrees. Then to finish the proof, it suffices to show that $\zeta$
is injective on the whole $ZHH_i(A)$.

Indeed, assume given some class $a \in ZHH_i(A)$ such that
$\zeta(a)=0$. Then we can choose a finite sum decomposition
$$
C(a) = \sum_j \alpha_j(\tau_i)
$$
for some maps $\alpha_j:A_{i,i} \to A$. Let
$$
a' = a - \sum_j \alpha_j(\wt{\tau}_i),
$$
where $\wt{\tau}_i \in ZHH_i(A_{i,i}) \subset HH_i(A_{i,i})$ is our
fixed lifting of $\tau_i$. Then $C(a')=0$, so that $a'$ lies in
$\xi(BHH_i(A)) \subset ZHH_i(A)$, and therefore $\zeta(a')$ lies in
$b'HH_i(A) = bHH_i(A)$. But since $\zeta(a) = 0$, this implies that
$$
\overline{a} = \sum_j \zeta(\alpha_j(\wt{\tau}_i)) = \zeta(a-a')
$$
lies in $bHH_\idot(A) = b'HH_i(A) \subset HH_i(A)$, so that if we
denote by $C_{cl}$ the classical Cartier map for the algebra $A$, we
have
$$
0 = C_{cl}(\overline{a}) = \sum_j
C_{cl}(\alpha_j(\zeta(\wt{\tau}_i))) = \sum_j
\alpha_j(C(\wt{\tau}_i)) = \sum_j \alpha_j(\tau_i) = C(a).
$$
Thus $a$ lies in $\xi(BHH_i(A)) \subset ZHH_i(A)$. Since we already
know that $\zeta$ is injective on $\xi(BHH_i(A))$, we have $a = 0$.
\endproof

\bigskip

\noindent
{\sc
Steklov Math Institute, Algebraic Geometry section\\
\mbox{}\hspace{30mm}and\\
Laboratory of Algebraic Geometry, NRU HSE\\
\mbox{}\hspace{30mm}and\\
Center for Geometry and Physics, IBS, Pohang, Rep. of Korea
}

\medskip

\noindent
{\em E-mail address\/}: {\tt kaledin@mi.ras.ru}

\end{document}